\documentclass[12pt]{amsart}
\usepackage{amsthm,amssymb}
\newcommand{\ol}{\overline}

\newcommand{\iy}{\infty}

\newcommand{\vf}{\varphi}

\newcommand{\al}{\alpha}

\newcommand{\dt}{\delta}

\newcommand{\sbs}{\subset}

\newcommand{\supp}{\operatorname{supp}}

\newtheorem{theorem}{\bf  Theorem}

\begin{document}
\begin{center}
    {\bf ON PARAMETRIZATION OF COMPACT WAVELET MATRICES}

\vskip+0.5cm
   {Lasha Ephremidze, \fbox{Gigla Janashia}\,, and Edem Lagvilava  }
\end{center}


\vskip+0.5cm

{\small {\bf Abstract.} We give an efficient complete
parametrization of wavelet matrices of rank $m$, genus $g+1$, and
degree $g$, which are naturally identified with corresponding
polynomial paraunitary  matrix-functions. The parametrization
depends on Wiener-Hopf factorization of unitary matrix-functions
with constant determinant given in the unit circle. This method
allows us to construct in real time the coefficients of wavelet
matrices from the above class.}

\vskip+1cm

A wavelet matrix ${\mathcal A}=(a_j^r)$ of rank $m$ consists of
$m$ rows of possibly infinite vectors
\begin{equation}
{\mathcal A}= \left(\begin{matrix} \cdots&
a^0_{-1}&a^0_{0}&a^0_{1}&a^0_{2}&\cdots\\[1mm]
\cdots& a^1_{-1}&a^1_{0}&a^1_{1}&a^1_{2}&\cdots\\
&\vdots&\vdots\\
\cdots& a^{m-1}_{-1}&a^{m-1}_{0}&a^{m-1}_{1}&a^{m-1}_{2}&\cdots\\
\end{matrix}\right),
\end{equation}
$a_j^r\in \mathbb{C}$, satisfying following two conditions.

(i) {\em Quadratic condition}:
\begin{equation}
 \sum_j a^r_{j+ml}\,\ol{a}^s_{j+mn}=m\,\dt^{rs}\dt_{nl};
\end{equation}

(ii) {\em Linear condition:}
\begin{equation}
\sum_{j=-\iy}^{\iy}a_j^r=m\dt^{r,0},
\end{equation}
where $\dt$ stands for the Kronecker symbol.

In this paper we assume that (1) is {\em compact}, i.e. only
finite number of its entries are different from $0$. Therefore,
the series in (2) and (3) are only formally infinite, and no
problem of convergence appears.

The quadratic condition (2) asserts that the rows of a wavelet
matrix ${\mathbf a}^r:=(a^r_j)_{j=\ol{-\iy,\iy}}$ have length
equal to $\sqrt{m}$ and that they are pairwise orthogonal when
shifted by an arbitrary multiple of $m$. The first row ${\mathbf
a}^0$ is called the {\em scaling vector} or {\em low-pass filter},
while remaining rows ${\mathbf a}^r$ are called the {\em wavelet
vectors} or {\em high-pass filters}. In signal processing
applications, the linear constraint (3) implies that a constant
signal emerges from the first subband of the {\em maltirate filter
bank} (1).

Associate to each wavelet matrix $\mathcal{A}$ the matrix function
$\mathbf{A}(z)$ as follows: let  $A_k$, $k\in\mathbb{Z}$, be
submatrices of $\mathcal{A}$ of size $m\!\times\!m$ defined by
$A_k=(a^r_{km+s})$, $0\leq s,r\leq m-1$, in other words, (1) is
expressed in terms of block matrices in the form
$$
{\mathcal A}=(\cdots, A_{-1},A_0,A_1,A_2,\cdots),
$$
and assume
\begin{equation}
{\mathbf A}(z)=\sum_{k=-\iy}^\iy A_kz^k\,.
\end{equation}
Obviously, there is one-to-one correspondence between matrices (1)
and formal series expansions (4), and,  for a compact  matrix, the
corresponding matrix-function is a Laurent polynomial.

It can be verified that the quadratic and the linear constraints
on ${\mathcal A}$ are equivalent, respectively, to the following
two conditions on ${\mathbf A}(z)$:
\begin{equation}
{\mathbf A}(z) {\mathbf A}^*(z^{-1})=mI,
\end{equation}
where ${\mathbf A}^*(z^{-1}):=\sum_{k=-\iy}^\iy A_k^*z^{-k}$ is
the adjoint of ${\mathbf A}(z)$, and
\begin{equation}
\sum_{j=1}^{m} {\mathbf A}_{ij}(1)=m\dt_{i,1},\;\;1\leq i\leq m,
\end{equation}
where ${\mathbf A(z)}=\big({\mathbf A}_{ij}(z)\big)_{i,j=1}^{m}$.
The condition (5) means that ${\mathbf A}$ is a {\em paraunitary
matrix-function}.

If $U$ is a  unitary matrix of size $m\!\times\!m$, $U\in
\mathcal{U}(m)$,  and ${\mathbf A}(z)$ satisfies (5), then
$U{\mathbf A}(z)$ satisfies (5) as well. Furthermore, for each
paraunitary matrix-function ${\mathbf A}(z)$, there exists and one
can explicitly construct a unitary matrix $U$ such that $U{\mathbf
A}(z)$ satisfies the linear condition (6) as well. If $U$ and $U'$
are two such matrices, then
$$
U'=\left(\begin{matrix}1&0\\0&V \end{matrix}\right) U,
$$
where $V\in\mathcal{U}(m-1)$. Thus the construction of paraunitary
matrix-functions are decisive for construction of wavelet
matrices.

It is said  that a wavelet matrix (1) has the {\em rank} $m$ and
the {\em genus} $g$, ${\mathcal A} \in WM(m,g;\mathbb{C})$, if the
corresponding matrix-function ${\mathbf A}(z)$ has a form
\begin{equation}
{\mathbf A}(z)=\sum_{k=0}^{g-1} A_kz^k.
\end{equation}
It can be easily shown (see [5, p. 58] that the determinant of a
paraunitary matrix-function ${\mathbf A}(z)$ is a monomial in $z$,
that is, there is a nonnegative integer $d$, called the {\em
degree}  of ${\mathbf A}(z)$, such that
$$
\det {\mathbf A}(z)=cz^d.
$$
Generically, (7) has degree $g-1$, although in specific
degenerated cases, it can be larger or smaller than $g-1$.

The relation between compact wavelet matrices  and  compactly
supported wavelet systems as orthonormal functions in
$L^2(\mathbb{R})$ is well-known (see [5], Ch. 5).

\smallskip

{\bf Theorem} ([1], [4], for rank 2; [5, pp. 87, 91], for rank
$m>2$): Let
$${\mathcal A}\in WM(m,g;\mathbb{C})$$
be a wavelet matrix and consider the functional difference
equation
\begin{equation}
\phi(x)=\sum_{k=0}^{mg-1}a_k^0\phi(mx-k)
\end{equation}
called the {\em scaling equation} associated with ${\mathcal A}$.
Then, there exists a unique $\phi\in L_2(\mathbb{R})$, called the
{\em scaling function}, which solves (8) and satisfies
$$\int\limits_\mathbb{R}\phi(x)\,dx=1\;\text{and }\;\;\supp\phi\sbs
\left[0,(g-1)\left(\frac m{m-1}\right)+1 \right].$$

Furthermore, if we define {\em wavelet functions} (associated with
${\mathcal A}$) by the formula
$$
\psi^r(x)=\sum_{k=0}^{mg-1}a_k^r\phi(mx-k),\;\;\;\;1\leq r<m,
$$
and consider the collection of functions
\begin{gather*}
\phi_{jk}(x)=m^{j/2}\phi(m^jx-k),\;\;j,k\in\mathbb{Z},\\
\psi^r_{jk}(x)=m^{j/2}\psi^r(m^jx-k),\;1\leq
r<m;\;j,k\in\mathbb{Z},
\end{gather*}
called the wavelet system ${\mathcal W}[{\mathcal A}]$ (associated
with wavelet matrix ${\mathcal A}$), then there exists an
$L_2$-convergent expansion for each $f\in L_2$:
\begin{equation}
f(x)=\sum_{k=-\iy}^\iy
c_{k}\phi_{0k}(x)+\sum_{r=1}^{m-1}\sum_{j=0}^\iy\sum_{k=-\iy}^\iy
c^r_{jk}\psi^r_{jk}(x),
\end{equation}
where the coefficients are given by
\begin{gather*}
c_{k}=\int_{-\iy}^\iy f(x)\ol{\phi_{0k}(x)}\,dx\,,\\
c^r_{jk}=\int_{-\iy}^\iy f(x)\ol{\psi^r_{jk}(x)}\,dx.
\end{gather*}

\smallskip

{\em Remark}: For most wavelet matrices ${\mathcal A}$, the
wavelet system ${\mathcal W}[{\mathcal A}]$ is a complete
orthonormal system and hence an orthonormal bases for
$L_2(\mathbb{R})$, which would imply the above theorem. However,
for some wavelet matrices, the system ${\mathcal W}[{\mathcal A}]$
is not orthonormal, and yet (9) is always true, which means that
${\mathcal W}[{\mathcal A}]$ is a {\em tight frame}.

\smallskip

Independently from the above mentioned connection between the
wavelet matrices and associated wavelet systems, the former can be
directly used in various discrete signal processing applications.
Namely, the following theorem is one of the key links between the
mathematical theory of wavelets and its practical applications.

\smallskip

{\bf Theorem} ({\bf wavelet matrix
 expansion}, [5, p. 80]): Let $$f:\mathbb{Z}\to\mathbb{C}$$ be an arbitrary
function (discrete signal) and let $${\mathcal A}=(a^r_k)\in
WM(m,g;\mathbb{C})$$ be a wavelet matrix of rank
 $m$ and genus $g$. Then $f$ has a unique  wavelet matrix
 expansion
 $$
f(n)=\sum_{r=0}^{m-1}\sum_{k=-\iy}^\iy c^r_ka^r_{mk+n}\,,
 $$
 where
 $$
c^r_k=\frac1m\sum_{n=-\iy}^\iy f(n)\ol{a}^r_{mk+n}\,.
 $$

 \smallskip

 The wavelet matrix expansion is locally finite; that is, for
 given $n$, only finitely many terms of the series are different
 from $0$.

 \smallskip

 From whatever said above it is evident the theoretical and practical
importance of
 deeper understanding an internal structure of paraunitary matrix-functions
which would
 allow to construct  efficiently a wide class of such matrices.
 So far, the only way of  classification of
 paraunitary matrix functions was via the following factorization
 theorem. This theorem resembles the factorization of polynomials
 of degree $d$ according to their $d$ roots and highest
 coefficients.

 For a unit column vector $v\in \mathbb{C}^m$, $v^*v=1$, let
 \begin{equation}
V(z):=I-vv^*+vv^*z.
\end{equation}
Obviously, (10) is a polynomial matrix function of order 1. It can
be shown that $V(z)$ is the paraunitary matrix-function of degree
1 (see [5, p. 59]) and it is called {\em primitive}.

\smallskip

{\bf Theorem} ({\bf Paraunitary Matrix Factorization}, [5, p.
60]): A paraunitary matrix-function (7) of  degree $d$, where
$A_{g-1}\not=0$ can be factorized as
$$
{\mathbf A}(z)=V_1(z)V_2(z)\cdots V_d(z)U
$$
where $V_j(z)$, $j=1,2,\ldots,d$, are primitive paraunitary
matrix-functions and $U$ is a (constant) unitary matrix.

\smallskip

We propose absolutely new way of parametrization of paraunitary
matrix-functions of rank $m$, genus $g+1$ and order $g$, which
depends on  Wiener-Hopf factorization of unitary matrix-functions
(with constant determinant) given on the unit circle in the
complex plane. Actually this method was developed in [3], [2] and
it  allows to construct efficiently matrix-functions of the above
type (consequently, to prepare the coefficients of the whole class
of compactly supported wavelets) in real time.

Let
\begin{gather*}
{\mathbf A}(z)=\begin{pmatrix} {\mathbf A}_{1,1}(z)& {\mathbf
A}_{1,2}(z)&
\cdots&{\mathbf A}_{1,m}(z)\\
{\mathbf A}_{2,1}(z)& {\mathbf A}_{2,2}(z)& \cdots&{\mathbf A}_{2,m}(z)\\
\vdots&\vdots&\vdots&\vdots\\{\mathbf A}_{m,1}(z)& {\mathbf
A}_{m,2}(z)& \cdots&{\mathbf A}_{m,m}(z)\end{pmatrix},\\
{\mathbf A}_{rj}(z)=\sum_{k=0}^{g}\al_k^{rj}\,z^k,\;\;\;1\leq
r,j\leq m,\\
{\mathbf A}(z) {\mathbf A}^*(z^{-1})=I,\\
\det {\mathbf A}(z)=cz^g\,,\;\;\;\;\;\;|c|=1.
\end{gather*}
We first convert the matrix-function ${\mathbf A}(z)$ into a
unitary (on the unit circle) matrix-function $U(z)$ by dividing
any row of ${\mathbf A}(z)$ (say, the last row, to be specific) by
$z^g$:
$$
U(z):=\begin{pmatrix} {\mathbf A}_{1,1}& {\mathbf A}_{1,2}&
\cdots&{\mathbf
A}_{1,m}\\
{\mathbf A}_{2,1}& {\mathbf A}_{2,2}& \cdots&{\mathbf A}_{2,m}\\
\vdots&\vdots&\vdots&\vdots\\
{\mathbf A}_{m-1,1}& {\mathbf A}_{m-1,2}& \cdots&{\mathbf
A}_{m-1,m}
\\z^{-g}{\mathbf A}_{m,1}&
z^{-g}{\mathbf A}_{m,2}& \cdots&z^{-g}{\mathbf
A}_{m,m}\end{pmatrix}.
$$
Then we have
\begin{gather*}
U(z)U^*(z)=I,\;\;\;\;\;\;\text{for }|z|=1, \\
\det U(z)=c,\;\;\;\;\;\;\;|c|=1,\\
 U_{rj}\in L^+_g\,,\;\;\;\;\;1\leq r<m;\;1\leq j\leq m,\\
U_{m,j}\in L^-_g\,,\;\;\;\;\;1\leq j\leq m,
\end{gather*}
where
\begin{gather*}
L^+_g=\big\{f: f(z)=\sum_{k=0}^{g} c_kz^k\big\},
\\
L^-_g=\big\{f: f(z)=\sum_{k=0}^{g} c_kz^{-k}\big\}.
\end{gather*}

Thus, the two theorems below give simple and transparent way of
one-to-one parametrization of paraunitary matrix-functions of rank
$m$, genus $g+1$ and degree $g$. To compare this with the above
presented factorization theorem from the simplicity point of view,
the proposed parametrization resembles the classification of
polynomials of degree $d$ according to their $d+1$ coefficients.

\begin{theorem} {\rm (see [2], p. 22)} For each $m\!\times\!m$ matrix-function
$F(z)$ of form
\begin{equation}
F=\begin{pmatrix}1&0&0&\cdots&0&0\\
          0&1&0&\cdots&0&0\\
           0&0&1&\cdots&0&0\\
           \vdots&\vdots&\vdots&\vdots&\vdots&\vdots\\
           0&0&0&\cdots&1&0\\
           \vf_{1}&\vf_{2}&\vf_{3}&\cdots&\vf_{m-1}&1
           \end{pmatrix},
\end{equation}
where
\begin{equation}
\vf_j\in L_g^-\,,\;\;\;\;1\leq j\leq m-1,
\end{equation}
there exists a unitary matrix-function $U(z)$ (unique up to a
constant unitary right multiplier) of form
\begin{equation}
U=\begin{pmatrix}u_{11}&u_{12}&\cdots&u_{1m}\\
                 u_{21}&u_{22}&\cdots&u_{2m}\\
           \vdots&\vdots&\vdots&\vdots\\
           u_{m-1,1}&u_{m-1,2}&\cdots&u_{m-1,m}\\
           u^*_{m1}&u^*_{m2}&\cdots&u^*_{mm}\\
           \end{pmatrix},
\end{equation}
 where
\begin{equation}
u_{kj}\in L_g^{+},\;\;\;1\leq k,j\leq m,
\end{equation}
with constant determinant, such that
\begin{equation}
F(z)U(z)\in L^+_g.
\end{equation}
\end{theorem}

\begin{theorem} {\rm (see [3])} For each unitary matrix-function $U(z)$, with
constant determinant,  of form $(13)$, $(14)$, there exists an
unique $F(z)$ of form $(11)$, $(12)$ such that $(15)$ holds.
\end{theorem}

Observe that if we denote by $F_-(z)$ the matrix function of type
(11) where each $\vf_j$ is replaced by $-\vf_j$,
$j=1,2,\ldots,m-1$, then $F_-(z)=(F(z))^{-1}$, so that the
equation
$$
U(z)=F_-(z)\cdot F(z)U(z)
$$
gives the rigt Wiener-Hopf factorization of $U(z)$.

\smallskip

In the end, we should mention that less than 1 sc computer time is
required to compute coefficients of matrix-function (13) in
Theorem 1 whenever coefficients of functions $\vf_j$,
$j=1,2,\ldots,m-1$, are selected in (11) for such large dimensions
as $m=30$ and $g=50$. This speed of calculations opens possibility
to choose the optimal wavelet matrix for specific problem, which
is the most important step in practical applications, by total
selection.

\end{document}